\documentclass{article}
\usepackage{amssymb}
\usepackage{amsfonts}
\usepackage{amsmath}
\usepackage[left=1.5in, right=1.5in, bottom=1.5in,top=1.5in]{geometry}

\input{tcilatex}
\begin{document}

\author{A. El- Gohary\thanks{%
Corresponding author. Permanent address: Department of Mathematics, College
of Science, Mansoura University, Mansoura 35516, Egypt. E-mail address:
elgohary0@yahoo.com (A. El-Gohary)}, M. El- Morshedy. \\
%EndAName
Mathematics Department, \\
Faculty of Science, Mansoura University, Mansoura,\\
Egypt.}
\title{Bivariate Exponentiated Modified Weibull Extension Distribution}
\date{}
\maketitle

\begin{abstract}
In this paper, \ we introduce a new bivariate distribution we called it
bivariate exponentiated modified Weibull extension distribution (BEMWE). The
model introduced here is of Marshall-Olkin type. The marginals of the new
bivariate distribution have exponentiated modified Weibull extension
distribution which proposed by Sarhan et al.(2013). The joint probability
density function and the joint cumulative distribution function are in
closed forms. Several properties of this distribution have been
discussed.The maximum likelihood estimators of the parameters are derived.
One real data set are analyzed using the new bivariate distribution, which
show that the new bivariate distribution can be used quite effectively in
fitting and analyzing real lifetime data.

\end{abstract}

\ \ \ \ \ \ \ \ \ \ \ \ 

\section{Introduction}

Recently, Sarhan et al. (2013) has defined a new four-parameter distribution
referred to as exponentiated modified Weibull extension (EMWE) distribution.
Sarhan et al. (2013) defined the (EMWE) distribution by exponentiating the
new modified Weibull extension (MWE) distribution which discussed by Xie et
al. (2002) as wos done for the exponentiated weibull (EW) distribution by
Mudholkar et al. (1995). They observed that exponential distribution,
generalized exponential distribution (1999), Gompertz distribution (1824),
generalized Gompertz (GG) distribution (2013), exponentiated Weibull (EW)
distribution (1995), Weibull extension model of Chen (2000), modified
Weibull extension (MWE) distribution (2002) and etc distribution can be
obtained as special cases of the (EMWE) distribution.

The objective of this paper is to provide a new bivariate distribution,
whose marginals are (EMWE) distributions which referred to as bivariate
exponentiated modified Weibull extension (BEMWE) distribution. It is
obtained using a method similar to that used to obtain Marshall-Olkin
bivariate exponential model Marshall and Olkin (1967).

The paper is organized as follows. Section 2 presents the shock model
yielding the (BEMWE) distribution. Also, the joint cumulative distribution
function, the joint probability density function, the marginal probability
density functions and the conditional probability density functions of
(BEMWE) distribution is derived in Section 2. In Section 3 sum reliability
studies are obtained. Section 4 presents the the marginal expectation of the
(BEMWE) distribution. Section 5 obtains the parameter estimation using MLE.
In section 6 a numerical result are obtained using real data. Finally, a
conclusion for the results is given in Section 7.

\section{\textbf{Bivariate exponentiated modified Weibull extension
distribution }}

In this section we introduce the BEMWE distribution using a method similar
to that which was used by Marshall and Olkin (1967) to define the Marshall
Olkin bivariate exponential (MOBE) distribution. We start with the joint
cumulative function of the proposed bivariate distribution and so used it to
derive the corresponding joint probability density function. Finally The
marginal probability density functions and conditional probability density
functions of this distribution are also derived. Let X be a random variable
has univariate EMWE distribution with parameters $\gamma ,\alpha ,\beta
,\lambda >0,$ then the corresponding cumulative distribution function (CDF)
is given by%
\begin{equation}
F(x)=\left[ 1-e^{-\lambda \alpha (e^{(x/\alpha )^{\beta }}-1)}\right]
^{\gamma },\text{ \ \ }x\geq 0,  \label{1}
\end{equation}%
and the probability density function (PDF) takes the following form%
\begin{equation}
f(x)={\small \gamma \lambda \beta e}^{(x/\alpha )^{\beta }}(\frac{{\small x}%
}{{\small \alpha }})^{\beta -1}e^{-\lambda \alpha (e^{(x/\alpha )^{\beta
}}-1)}\left[ 1-e^{-\lambda \alpha (e^{(x/\alpha )^{\beta }}-1)}\right]
^{\gamma -1},\text{ \ \ }x\geq 0.  \label{2}
\end{equation}

\subsection{Joint cumulative distribution function}

Suppose that $U_{i}$ $(i=1,2,3)$ are three independent random variables such
that $U_{i}$ $\sim $EMWE ($\gamma _{i},\alpha ,\beta ,\lambda $). Define $%
X_{1}=max\{U_{1},U_{3}\}$ and $X_{2}=max\{U_{2},U_{3}\}$. Then we say that
the bivariate vector $(X_{1},X_{2})$ has a bivariate exponentiated modified
Weibull extension distribution, with parameters $(\gamma _{1},\gamma
_{2},\gamma _{3},\alpha ,\beta ,\lambda )$ and we denote it by BEMWE$(\gamma
_{1},\gamma _{2},\gamma _{3},\alpha ,\beta ,\lambda ).$ The following
interpretation can be provided for the BEMWE model.

\textit{Shock model:} Assum thate there exists a three independent sources
of shocks. Suppose these shocks are affecting a system with two components.
It is assumed that the shock from source 1 reaches the system and destroys
component 1 immediately, the shock from source 2 reaches the system and
destroys component 2 immediately, while if the shock from source 3 hits the
system it destroys both the components immediately. Let $U_{i}$ denote the
inter-arrival times, between the shocks in source $i,$ $i=1,2,3,$ which
follow the distribution EMWE. If $X_{1},$ $X_{2}$ denote the survival times
of the components, then the bivariate vector $(X_{1},X_{2})$ follows the
BEMWE model.

We now study the joint cumulative distribution function of the bivariate
random vector $(X_{1},X_{2})$ in the following lemma.

\paragraph{\textbf{Lemma 2.1. }}

The joint CDF of $(X_{1},X_{2})$ is%
\begin{equation}
F_{BEMWE}(x_{1},x_{2})=\left[ 1-e^{-\lambda \alpha (e^{(x_{1}/\alpha
)^{\beta }}-1)}\right] ^{\gamma _{1}}\left[ 1-e^{-\lambda \alpha
(e^{(x_{2}/\alpha )^{\beta }}-1)}\right] ^{\gamma _{2}}\left[ 1-e^{-\lambda
\alpha (e^{(z/\alpha )^{\beta }}-1)}\right] ^{\gamma _{3}},  \label{3}
\end{equation}%
where $z=min\left( x_{1},x_{2}\right) .$

\paragraph{\textbf{proof: }\ }

Since the joint CDF of the random variables $X_{1}$ and $X_{2}$ is defined as%
\begin{eqnarray*}
F(x_{1},x_{2}) &=&P\left( X_{1}\leq x_{1},X_{2}\leq x_{2}\right) \\
&=&P\left( max\{U_{1},U_{3}\}\leq x_{1},max\{U_{2},U_{3}\}\leq x_{2}\right)
\\
&=&P\left( U_{1}\leq x_{1},U_{2}\leq x_{2},U_{3}\leq min\left(
x_{1},x_{2}\right) \right) .
\end{eqnarray*}%
As the random variables $U_{i}$ $(i=1,2,3)$ are mutually independent, we
directly obtain%
\begin{equation*}
F_{BEMWE}(x_{1},x_{2})=P\left( U_{1}\leq x_{1}\right) P\left( U_{2}\leq
x_{2}\right) P\left( U_{3}\leq min\left( x_{1},x_{2}\right) \right) \text{ \
\ \ \ \ \ \ \ \ \ \ \ \ \ \ \ \ \ \ \ \ \ \ \ \ \ \ \ \ }
\end{equation*}%
\begin{equation}
=F_{EMWE}(x_{1};\gamma _{1},\alpha ,\beta ,\lambda )F_{EMWE}(x_{2};\gamma
_{2},\alpha ,\beta ,\lambda )F_{EMWE}(z;\gamma _{3},\alpha ,\beta ,\lambda )
\label{4}
\end{equation}%
Substituting from (1) into (4), we obtain (3), which completes the proof of
the lemma.2.1.

\subsection{Joint probability density function}

The following theorem gives the joint PDF of the $X_{1}$ and $X_{2}$ which
is the joint PDF of BEMWE $(\gamma _{1},\gamma _{2},\gamma _{3},\alpha
,\beta ,\lambda )$.

\paragraph{\textbf{Theorem 2.1.}}

If the joint CDF of $X_{1}$ and $X_{2}$ is as in (3), then the joint PDF of $%
X_{1}$ and $X_{2}$ takes the form%
\begin{equation}
f_{BEMWE}(x_{1},x_{2})=\left\{ 
\begin{array}{l}
f_{1}(x_{1},x_{2})\ \ \ \ \text{if \ }x_{1}<x_{2} \\ 
f_{2}(x_{1},x_{2})\ \ \ \ \text{if \ }x_{2}<x_{1} \\ 
f_{0}(x,x)\ \ \ \ \ \ \ \text{if\ \ }x_{1}=x_{2}=x%
\end{array}%
,\right.  \label{5........}
\end{equation}%
where%
\begin{eqnarray}
f_{1}(x_{1},x_{2}) &=&f_{EMWE}(x_{2};\gamma _{2},\alpha ,\beta ,\lambda
)f_{EMWE}(x_{1};\gamma _{1}+\gamma _{3},\alpha ,\beta ,\lambda )  \notag \\
&=&{\small \gamma }_{2}\left( \gamma _{1}+\gamma _{3}\right) {\small \lambda 
}^{2}{\small \beta }^{2}{\small e}^{(x_{2}/\alpha )^{\beta }}(\frac{{\small x%
}_{2}}{{\small \alpha }})^{\beta -1}e^{-\lambda \alpha (e^{(x_{2}/\alpha
)^{\beta }}-1)}\left[ 1-e^{-\lambda \alpha (e^{(x_{2}/\alpha )^{\beta }}-1)}%
\right] ^{\gamma _{2}-1}  \notag \\
&&{\small \times e}^{(x_{1}/\alpha )^{\beta }}(\text{$\frac{{\small x}_{1}}{%
{\small \alpha }})$}^{{\small \beta }-1}e^{-\lambda \alpha (e^{(x_{1}/\alpha
)^{\beta }}-1)}\left[ {\tiny 1-}e^{-\lambda \alpha (e^{(x_{1}/\alpha
)^{\beta }}-1)}\right] ^{\gamma _{1}+\gamma _{3}-1},  \label{6..}
\end{eqnarray}%
\begin{eqnarray}
f_{2}(x_{1},x_{2}) &=&f_{EMWE}(x_{1};\gamma _{1},\alpha ,\beta ,\lambda
)f_{EMWE}(x_{2};\gamma _{2}+\gamma _{3},\alpha ,\beta ,\lambda )  \notag \\
&=&{\small \gamma }_{1}\left( \gamma _{2}+\gamma _{3}\right) {\small \lambda 
}^{2}{\small \beta }^{2}{\small e}^{(x_{1}/\alpha )^{\beta }}(\frac{x_{1}}{%
\alpha })^{\beta -1}e^{-\lambda \alpha (e^{(x_{1}/\alpha )^{\beta }}-1)}%
\left[ 1-e^{-\lambda \alpha (e^{(x_{1}/\alpha )^{\beta }}-1)}\right]
^{\gamma _{1}-1}  \notag \\
&&{\tiny \times e}^{(x_{2}/\alpha )^{\beta }}(\frac{{\tiny x}_{_{2}}}{{\tiny %
\alpha }})^{{\small \beta -1}}e^{-\lambda \alpha (e^{(x_{2}/\alpha )^{\beta
}}-1)}\left[ 1-e^{-\lambda \alpha (e^{(x_{2}/\alpha )^{\beta }}-1)}\right]
^{\gamma _{2}+\gamma _{3}-1}  \label{7...}
\end{eqnarray}%
and%
\begin{eqnarray}
f_{3}(x,x) &=&\frac{\gamma _{3}}{\gamma _{1}+\gamma _{2}+\gamma _{3}}%
f_{EMWE}(x_{2};\gamma _{1}+\gamma _{2}+\gamma _{3},\alpha ,\beta ,\lambda )%
\text{ \ }  \notag \\
&=&\gamma _{3}\lambda \beta e^{(x/\alpha )^{\beta }}(\frac{x}{\alpha }%
)^{\beta -1}e^{-\lambda \alpha (e^{(x/\alpha )^{\beta }}-1)}\left[
1-e^{-\lambda \alpha (e^{(x/\alpha )^{\beta }}-1)}\right] ^{\gamma
_{1}+\gamma _{2}+\gamma _{3}-1}.\text{ }  \label{8..}
\end{eqnarray}

\paragraph{Proof}

Let us first assume that $x_{1}<x_{2}.$ Then, the expression for $%
f_{1}(x_{1},x_{2})$ can be simply obtained by differentiating the joint CDF $%
F_{BEMWE}(x_{1},x_{2})$ given in (3) with respect to $x_{1}$ and $x_{2}.$
Simillary, we find the expression of $f_{2}(x_{1},x_{2})$ when $x_{2}<x_{1}.$
But $f_{3}(x,x)$ can not be derived in a similar method. For this reason, we
use the following identity to derive $f_{3}(x,x).$ 
\begin{equation}
\dint\limits_{0}^{\infty
}\dint\limits_{0}^{x_{2}}f_{1}(x_{1},x_{2})dx_{1}dx_{2}+\dint\limits_{0}^{%
\infty
}\dint\limits_{0}^{x_{1}}f_{2}(x_{1},x_{2})dx_{2}dx_{1}+\dint\limits_{0}^{%
\infty }f_{3}(x,x)dx=1  \label{9..}
\end{equation}%
Let%
\begin{equation*}
I_{1}=\dint\limits_{0}^{\infty
}\dint\limits_{0}^{x_{2}}f_{1}(x_{1},x_{2})dx_{1}dx_{2}\text{ \ \ \ and \ }%
I_{2}=\dint\limits_{0}^{\infty
}\dint\limits_{0}^{x_{1}}f_{2}(x_{1},x_{2})dx_{2}dx_{1}
\end{equation*}%
One can find that%
\begin{equation}
I_{1}=\dint\limits_{0}^{\infty }\gamma _{2}{\small \lambda \beta e}%
^{(x_{2}/\alpha )^{\beta }}(\frac{{\small x}_{2}}{{\small \alpha }})^{\beta
-1}e^{-\lambda \alpha (e^{(x_{2}/\alpha )^{\beta }}-1)}\left[ 1-e^{-\lambda
\alpha (e^{(x_{2}/\alpha )^{\beta }}-1)}\right] ^{\gamma _{1}+\gamma
_{2}+\gamma _{3}-1}dx_{2}  \label{10..}
\end{equation}%
and%
\begin{equation}
I_{2}=\dint\limits_{0}^{\infty }\gamma _{1}{\small \lambda \beta e}%
^{(x_{1}/\alpha )^{\beta }}(\frac{{\small x}_{1}}{{\small \alpha }})^{\beta
-1}e^{-\lambda \alpha (e^{(x_{1}/\alpha )^{\beta }}-1)}\left[ 1-e^{-\lambda
\alpha (e^{(x_{1}/\alpha )^{\beta }}-1)}\right] ^{\gamma _{1}+\gamma
_{2}+\gamma _{3}-1}dx_{1}  \label{11...}
\end{equation}%
Substituting from (10) and (11) into (9) we obtain%
\begin{equation*}
\dint\limits_{0}^{\infty }f_{3}(x,x)dx=1-I_{1}-I_{2}\text{ \ \ \ \ \ \ \ \ \
\ \ \ \ \ \ \ \ \ \ \ \ \ \ \ \ \ \ \ \ \ \ \ \ \ \ \ \ \ \ \ \ \ \ \ \ \ \
\ \ \ \ \ \ \ \ \ \ \ \ \ \ \ \ \ \ \ \ \ \ \ \ \ \ \ \ \ \ \ \ \ \ }
\end{equation*}%
\begin{eqnarray*}
&=&\dint\limits_{0}^{\infty }\left( {\small \gamma }_{1}{\tiny +}{\small %
\gamma }_{2}{\tiny +}{\small \gamma }_{3}\right) {\small \lambda \beta e}%
^{\left( x/\alpha \right) ^{\beta }}{\small (}\frac{{\small x}}{{\small %
\alpha }}{\small )}^{{\small \beta -1}}e^{-\lambda \alpha (e^{(x/\alpha
)^{\beta }}-1)}\left[ {\small 1-}e^{-\lambda \alpha (e^{(x/\alpha )^{\beta
}}-1)}\right] ^{\gamma _{1}+\gamma _{2}+\gamma _{3}-1}{\small dx} \\
&&-\dint\limits_{0}^{\infty }\gamma _{2}\lambda \beta e^{\left( x/\alpha
\right) ^{\beta }}(\frac{{\small x}}{{\small \alpha }})^{\beta
-1}e^{-\lambda \alpha (e^{(x/\alpha )^{\beta }}-1)}\left[ 1-e^{-\lambda
\alpha (e^{(x/\alpha )^{\beta }}-1)}\right] ^{\gamma _{1}+\gamma _{2}+\gamma
_{3}-1}dx \\
&&-\dint\limits_{0}^{\infty }\gamma _{1}\lambda \beta e^{\left( x/\alpha
\right) ^{\beta }}(\frac{{\small x}}{{\small \alpha }})^{\beta
-1}e^{-\lambda \alpha (e^{(x/\alpha )^{\beta }}-1)}\left[ 1-e^{-\lambda
\alpha (e^{(x/\alpha )^{\beta }}-1)}\right] ^{\gamma _{1}+\gamma _{2}+\gamma
_{3}-1}dx.
\end{eqnarray*}%
Thus,%
\begin{equation*}
f_{3}(x,x)=\gamma _{3}\lambda \beta e^{(x/\alpha )^{\beta }}(\frac{x}{\alpha 
})^{\beta -1}e^{-\lambda \alpha (e^{(x/\alpha )^{\beta }}-1)}\left[
1-e^{-\lambda \alpha (e^{(x/\alpha )^{\beta }}-1)}\right] ^{\gamma
_{1}+\gamma _{2}+\gamma _{3}-1},
\end{equation*}

\subsection{Marginal probability density functions}

The following theorem gives the marginal probability density functions of $%
X_{1}$ and $X_{2}$.

\paragraph{Theorem 2.2.}

The marginal probability density functions of \ $X_{i}\ $\ ,$(i=1,2)$ is
given by%
\begin{eqnarray}
f_{X_{i}}(x_{i}) &=&f_{EMWE}(x_{i};\gamma _{i}+\gamma _{3},\alpha ,\beta
,\lambda ),\text{ \ \ }x_{i}>0,\text{ }i=1,2  \notag \\
&=&\left( \gamma _{i}+\gamma _{3}\right) \lambda \beta e^{(x_{i}/\alpha
)^{\beta }}(\frac{x_{i}}{\alpha })^{\beta -1}e^{-\lambda \alpha
(e^{(x_{i}/\alpha )^{\beta }}-1)}\left[ 1-e^{-\lambda \alpha
(e^{(x_{i}/\alpha )^{\beta }}-1)}\right] ^{\left( \gamma _{i}+\gamma
_{3}\right) -1}.  \label{12..}
\end{eqnarray}

\paragraph{Proof:}

The marginal cumulative distribution function for $X_{i}$ is 
\begin{equation*}
F(x_{i})=P\left( X_{i}\leq x_{i}\right) =P\left( max\{U_{i},U_{3}\}\leq
x_{i}\right) =P\left( U_{i}\leq x_{i},U_{3}\leq x_{i}\right) .
\end{equation*}%
As the random variables $U_{i}$ $(i=1,2)$ and $U_{3}$ are mutually
independent, we directly obtain%
\begin{eqnarray}
F(x_{i}) &=&P\left( U_{i}\leq x_{i})P(U_{3}\leq x_{i}\right)  \notag \\
&=&\left[ 1-e^{-\lambda \alpha (e^{(x_{i}/\alpha )^{\beta }}-1)}\right]
^{\gamma _{i}}\left[ 1-e^{-\lambda \alpha (e^{(x_{i}/\alpha )^{\beta }}-1)}%
\right] ^{\gamma _{3}}  \notag \\
&=&\left[ 1-e^{-\lambda \alpha (e^{(x_{i}/\alpha )^{\beta }}-1)}\right]
^{\gamma _{i}+\gamma _{3}}=F_{EMWE}(x_{i};\gamma _{i}+\gamma _{3},\alpha
,\beta ,\lambda ).  \label{13....}
\end{eqnarray}%
From which we readily derive the pdf of $X_{i}$ , $f(x_{i})=\frac{\partial }{%
\partial x_{i}}F(x_{i}),$ as in(12).

\subsection{Conditional probability density functions}

The following theorem gives the marginal probability density functions of $%
(X_{1},X_{2})$.

\paragraph{Theorem 2.3.}

The conditional probability density function of $X_{i}$ given $X_{j}=x_{j}$ $%
,$ $(i,j=1,2,i\neq j)$ is given by%
\begin{equation*}
f_{X_{i}\mid X_{j}}(x_{i}\mid x_{j})=\left\{ 
\begin{array}{l}
f_{X_{i}\mid X_{j}}^{(1)}(x_{i}\mid x_{j})\ \ \ \ \text{if }\ 0<x_{i}<x_{j}
\\ 
f_{X_{i}\mid X_{j}}^{(2)}(x_{i}\mid x_{j})\ \ \ \ \text{if \ }0<x_{j}<x_{i}
\\ 
f_{X_{i}\mid X_{j}}^{(3)}(x_{i}\mid x_{j})\ \ \ \ \text{if\ \ }x_{i}=x_{j}>0%
\end{array}%
.\right.
\end{equation*}%
where%
\begin{eqnarray*}
f_{X_{i}\mid X_{j}}^{(1)}(x_{i} &\mid &x_{j})=\frac{\gamma _{j}(\gamma
_{i}+\gamma _{3}){\small \lambda \beta e}^{(x_{i}/\alpha )^{\beta }}(\frac{%
{\small x}_{i}}{{\small \alpha }})^{\beta -1}e^{-\lambda \alpha
(e^{(x_{i}/\alpha )^{\beta }}-1)}\left[ 1-e^{-\lambda \alpha
(e^{(x_{i}/\alpha )^{\beta }}-1)}\right] ^{\gamma _{i}+\gamma _{3}-1}}{%
(\gamma _{j}+\gamma _{3})\left[ 1-e^{-\lambda \alpha (e^{(x_{j}/\alpha
)^{\beta }}-1)}\right] ^{\gamma _{i}+\gamma _{3}-1}}, \\
f_{X_{i}\mid X_{j}}^{(2)}(x_{i} &\mid &x_{j})=\gamma _{i}{\small \lambda
\beta e}^{(x_{i}/\alpha )^{\beta }}(\frac{{\small x}_{i}}{{\small \alpha }}%
)^{\beta -1}e^{-\lambda \alpha (e^{(x_{i}/\alpha )^{\beta }}-1)}\left[
1-e^{-\lambda \alpha (e^{(x_{i}/\alpha )^{\beta }}-1)}\right] ^{\gamma
_{i}-1}
\end{eqnarray*}%
and%
\begin{equation*}
f_{X_{i}\mid X_{j}}^{(3)}(x_{i}\mid x_{j})\ =\frac{\gamma _{3}}{\gamma
_{i}+\gamma _{3}}\left[ 1-e^{-\lambda \alpha (e^{(x_{i}/\alpha )^{\beta
}}-1)}\right] ^{\gamma _{i}}.
\end{equation*}

\paragraph{Proof:}

The proof follows immediately by substituting the joint probability density
function of $(X_{1},X_{2})$ given in (6), (7) and (8) and the marginal
probability density function of given in (12), using the relation%
\begin{equation*}
f_{X_{i}\mid X_{j}}(x_{i}\mid x_{j})\ =\frac{f_{X_{i},X_{j}}(x_{i},x_{j})}{%
f_{X_{i}}(x_{i})\ },\text{ \ }(i=1,2).
\end{equation*}

\section{Reliability studies}

In this section, we present the joint survival function of $(X_{1},X_{2})$,
the CDF of the random variable $Y=max\{X_{1},X_{2}\}$ and the CDF of the
random variable $W=\min \{X_{1},X_{2}\}.$

\subsection{Joint survival function}

In this subsection, we derive the joint survival function of $(X_{1},X_{2})$
in a compact form.

\paragraph{Theorem 3.1.}

The joint survival function of $(X_{1},X_{2})$ is given by%
\begin{equation}
S_{X_{1},X_{2}}(x_{1},x_{2})=\left\{ 
\begin{array}{l}
S_{1}(x_{1},x_{2})\ \ \ \ \text{if \ }x_{1}<x_{2} \\ 
S_{2}(x_{1},x_{2})\ \ \ \ \text{if \ }x_{2}<x_{1} \\ 
S_{0}(x,x)\ \ \ \ \ \ \ \text{if\ \ }x_{1}=x_{2}=x%
\end{array}%
,\right.  \label{14...}
\end{equation}%
where%
\begin{eqnarray*}
S_{1}(x_{1},x_{2}) &=&1-\left[ 1-e^{-\lambda \alpha (e^{(x_{2}/\alpha
)^{\beta }}-1)}\right] ^{\gamma _{2}+\gamma _{3}}-\left[ 1-e^{-\lambda
\alpha (e^{(x_{1}/\alpha )^{\beta }}-1)}\right] ^{\gamma _{1}+\gamma _{3}} \\
&&\times \left( 1-\left[ 1-e^{-\lambda \alpha (e^{(x_{2}/\alpha )^{\beta
}}-1)}\right] ^{\gamma _{2}}\right) ,
\end{eqnarray*}%
\begin{eqnarray*}
S_{2}(x_{1},x_{2}) &=&1-\left[ 1-e^{-\lambda \alpha (e^{(x_{1}/\alpha
)^{\beta }}-1)}\right] ^{\gamma _{1}+\gamma _{3}}-\left[ 1-e^{-\lambda
\alpha (e^{(x_{2}/\alpha )^{\beta }}-1)}\right] ^{\gamma _{2}+\gamma _{3}} \\
&&\times \left( 1-\left[ 1-e^{-\lambda \alpha (e^{(x_{1}/\alpha )^{\beta
}}-1)}\right] ^{\gamma _{1}}\right)
\end{eqnarray*}%
and%
\begin{eqnarray*}
S_{0}(x,x)\ &=&1-\left[ 1-e^{\lambda \alpha (1-e^{(x/\alpha )^{\beta }})}%
\right] ^{\gamma _{3}}\times \\
&&\left( \left[ 1-e^{-\lambda \alpha (e^{(x/\alpha )^{\beta }}-1)}\right]
^{\gamma _{1}}+\left[ 1-e^{-\lambda \alpha (e^{(x/\alpha )^{\beta }}-1)}%
\right] ^{\gamma _{2}}-\left[ 1-e^{-\lambda \alpha (e^{(x/\alpha )^{\beta
}}-1)}\right] ^{\gamma _{1}+\gamma _{2}}\right) .
\end{eqnarray*}

\paragraph{Proof:}

The joint survival function of $(X_{1},X_{2})$ can be obtained from the
following relation%
\begin{equation}
S_{X_{1},X_{2}}(x_{1},x_{2})=1-F_{X_{1}}(x_{1})-F_{X_{2}}(x_{2})+F_{X_{1},X_{2}}(x_{1},x_{2}).
\label{15...}
\end{equation}%
Substituting from (3) and (13) in (15), we get%
\begin{eqnarray}
S_{X_{1},X_{2}}(x_{1},x_{2}) &=&1-\left[ 1-e^{-\lambda \alpha
(e^{(x_{1}/\alpha )^{\beta }}-1)}\right] ^{\gamma _{1}+\gamma _{3}}-\left[
1-e^{-\lambda \alpha (e^{(x_{2}/\alpha )^{\beta }}-1)}\right] ^{\gamma
_{2}+\gamma _{3}}+  \notag \\
&&\left[ 1-e^{-\lambda \alpha (e^{(x_{1}/\alpha )^{\beta }}-1)}\right]
^{\gamma _{1}}\left[ 1-e^{-\lambda \alpha (e^{(x_{2}/\alpha )^{\beta }}-1)}%
\right] ^{\gamma _{2}}\left[ 1-e^{-\lambda \alpha (e^{(z/\alpha )^{\beta
}}-1)}\right] ^{\gamma _{3}},  \label{16....}
\end{eqnarray}%
where $z=min\left( x_{1},x_{2}\right) .$ From (16) we can be\ obtained
simply the expressions of $S_{1}(x_{1},x_{2}),$ $S_{2}(x_{1},x_{2})$ and $%
S_{0}(x_{1},x_{2})$ for $x_{1}<x_{2}$ , $x_{2}<x_{1}$ and $x_{1}=x_{2}=x$
respectively, which completes the proof.

\paragraph{Comment 3.1.}

Basu (1971) defined the bivariate failure rate function $h(x_{1},x_{2})$ for
the random vector $(X_{1},X_{2})$ as the following relation%
\begin{equation}
h_{X_{1},X_{2}}(x_{1},x_{2})=\frac{f_{X_{1},X_{2}}(x_{1},x_{2})}{%
S_{X_{1},X_{2}}(x_{1},x_{2})}.  \label{17....}
\end{equation}%
We can obtained the bivariate failure rate function $h(x_{1},x_{2})$ for the
random vector $(X_{1},X_{2})$ by substituting from (5) and (14) in (17).

\paragraph{\textbf{Lemma 3.1. }}

The CDF of the random variable $Y=max\{X_{1},X_{2}\}$ is given as%
\begin{equation}
F_{Y}(y)=\left[ 1-e^{-\lambda \alpha (e^{(y/\alpha )^{\beta }}-1)}\right]
^{\gamma _{1}+\gamma _{2}+\gamma _{3}}.  \label{18....}
\end{equation}%
\ \ 

\paragraph{Proof:}

Since%
\begin{eqnarray*}
F_{Y}(y) &=&P\left( Y\leq y\right) =P\left( max\{X_{1},X_{2}\}\leq y\right)
=P\left( X_{1}\leq y,X_{2}\leq y\right) \\
&=&P\left( max\{U_{1},U_{3}\}\leq y,max\{U_{2},U_{3}\}\leq y\right) =P\left(
U_{1}\leq y,U_{2}\leq y,U_{3}\leq y\right) ,
\end{eqnarray*}%
where the random variables $U_{i}$ $(i=1,2,3)$ are mutually independent, we
directly obtain%
\begin{eqnarray}
F_{Y}(y) &=&P\left( U_{1}\leq y)P(U_{2}\leq y)P(U_{3}\leq y\right)  \notag \\
&=&F_{EMWE}(y;\gamma _{1},\alpha ,\beta ,\lambda )F_{EMWE}(y;\gamma
_{2},\alpha ,\beta ,\lambda )F_{EMWE}(y;\gamma _{3},\alpha ,\beta ,\lambda ).
\label{19....}
\end{eqnarray}%
Substituting from (1) in (19), we get (18) which completes the proof of the
lemma.3.1.

\paragraph{Comment 3.2.}

From lemma 3.1. we can say that, if $X_{1}$ and $X_{2}$ are independent EMWE
random variables then $\ max\{X_{1},X_{2}\}$ is also EMWE random variable.

\paragraph{\textbf{Lemma 3.2. }}

The CDF of the random variable $W=\min \{X_{1},X_{2}\}$ is given as%
\begin{eqnarray}
F_{W}(w) &=&\left[ 1-e^{-\lambda \alpha (e^{(w/\alpha )^{\beta }}-1)}\right]
^{\gamma _{1}}+\left[ 1-e^{-\lambda \alpha (e^{(w/\alpha )^{\beta }}-1)}%
\right] ^{\gamma _{2}}  \notag \\
&&-\left[ 1-e^{-\lambda \alpha (e^{(w/\alpha )^{\beta }}-1)}\right] ^{\gamma
_{1}+\gamma _{2}+\gamma _{3}}.  \label{20.....}
\end{eqnarray}%
\ \ 

\paragraph{Proof:}

Since%
\begin{eqnarray}
F_{W}(w) &=&P\left( W\leq w\right) =P\left( \min \{X_{1},X_{2}\}\leq
w\right) =1-P\left( \min \{X_{1},X_{2}\}>w\right)  \notag \\
&=&1-P\left( X_{1}>w,X_{2}>w\right) =1-S(w,w)  \label{21.}
\end{eqnarray}%
Substituting from (14) in (21), we get%
\begin{equation}
F_{W}(w)=F_{X_{1}}(w)+F_{X_{2}}(w)-F_{X_{1},X_{2}}(w,w).  \label{22...}
\end{equation}%
Substituting from (3) and (13) in (22), we get (20) which completes the
proof of the lemma 3.2.

\section{The marginal expectation}

In this section, we derive the marginal expectation of $X_{i}\ (i=1,2).$ The
following theorem gives the $rth$ moments of $X_{i}\ (i=1,2)$ as infinite
series expansion.

\paragraph{Theorem 3.1.}

The $rth$ moment of \ $X_{i}\ (i=1,2)$ is given by:%
\begin{equation}
E(X_{i}^{r}\ )=\frac{\left( \gamma _{i}+\gamma _{3}\right) \lambda }{\alpha
^{\beta -1}}\dsum\limits_{j=0}^{\infty }\dsum\limits_{k=0}^{\infty }\binom{%
\left( \gamma _{i}+\gamma _{3}\right) -1}{j}\frac{\left( -1\right)
^{j+k}\lambda ^{k}\alpha ^{k+\beta +r}(j+1)^{k}}{(k+1)^{\beta +r}k!}%
e^{\lambda \alpha (j+1)}\Gamma (\frac{r}{\beta }+1).  \label{23...}
\end{equation}%
\ \ 

\paragraph{Proof:}

We will start with the known definition of the $rth$ moment of the random
variables $X_{i}$ with pdf $f(x_{i})$ given by%
\begin{equation*}
E(X_{i}^{r}\ )=\dint\limits_{0}^{\infty }x_{i}^{r}f_{X_{i}}(x_{i})dx_{i}.
\end{equation*}%
Substituting for $f_{X_{i}}(x_{i})$ from (12), we get%
\begin{equation}
E(X_{i}^{r}\ )=\frac{\left( \gamma _{i}+\gamma _{3}\right) \lambda \beta }{%
\alpha ^{\beta -1}}\dint\limits_{0}^{\infty }x_{i}^{r+\beta
-1}e^{(x_{i}/\alpha )^{\beta }}e^{-\lambda \alpha (e^{(x_{i}/\alpha )^{\beta
}}-1)}\left[ 1-e^{-\lambda \alpha (e^{(x_{i}/\alpha )^{\beta }}-1)}\right]
^{\left( \gamma _{i}+\gamma _{3}\right) -1}dx_{i}.  \label{24.}
\end{equation}%
Since $0<e^{-\lambda \alpha (e^{(x_{i}/\alpha )^{\beta }}-1)}<1$ for $x>0$,
then by using the binomial series expansion of $\left[ 1-e^{-\lambda \alpha
(e^{(x_{i}/\alpha )^{\beta }}-1)}\right] ^{\left( \gamma _{i}+\gamma
_{3}\right) -1}$given by%
\begin{equation}
\left[ 1-e^{-\lambda \alpha (e^{(x_{i}/\alpha )^{\beta }}-1)}\right]
^{\left( \gamma _{i}+\gamma _{3}\right) -1}=\dsum\limits_{j=0}^{\infty }%
\binom{\left( \gamma _{i}+\gamma _{3}\right) -1}{j}\left( -1\right)
^{j}e^{-j\lambda \alpha (e^{(x_{i}/\alpha )^{\beta }}-1)}.  \label{25..}
\end{equation}%
Substituting from (25) into (24), we get%
\begin{equation*}
E(X_{i}^{r}\ )=\frac{\left( \gamma _{i}+\gamma _{3}\right) \lambda \beta }{%
\alpha ^{\beta -1}}\dsum\limits_{j=0}^{\infty }\binom{\left( \gamma
_{i}+\gamma _{3}\right) -1}{j}\left( -1\right) ^{j}e^{\lambda \alpha
(j+1)}\dint\limits_{0}^{\infty }x_{i}^{r+\beta -1}e^{(x_{i}/\alpha )^{\beta
}}e^{-\lambda \alpha (j+1)e^{(x_{i}/\alpha )^{\beta }}}dx_{i}.
\end{equation*}%
Using the series expansion of \ $e^{-\lambda \alpha (j+1)e^{(x_{i}/\alpha
)^{\beta }}}$, one gets%
\begin{eqnarray*}
E(X_{i}^{r}\ ) &=&\frac{\left( \gamma _{i}+\gamma _{3}\right) \lambda \beta 
}{\alpha ^{\beta -1}}\dsum\limits_{j=0}^{\infty }\dsum\limits_{k=0}^{\infty }%
\binom{\left( \gamma _{i}+\gamma _{3}\right) -1}{j}\frac{\left( -1\right)
^{j+k}\lambda ^{k}\alpha ^{k}(j+1)^{k}}{k!}e^{\lambda \alpha (j+1)} \\
&&\times \dint\limits_{0}^{\infty }x_{i}^{r+\beta -1}e^{(k+1)(x_{i}/\alpha
)^{\beta }}dx_{i}.
\end{eqnarray*}%
Let $y=(k+1)(x_{i}/\alpha )^{\beta }$ in the above integral, then we can get%
\begin{equation}
E(X_{i}^{r}\ )=\frac{\left( \gamma _{i}+\gamma _{3}\right) \lambda }{\alpha
^{\beta -1}}\dsum\limits_{j=0}^{\infty }\dsum\limits_{k=0}^{\infty }\binom{%
\left( \gamma _{i}+\gamma _{3}\right) -1}{j}\frac{\left( -1\right)
^{j+k}\lambda ^{k}\alpha ^{k+\beta +r}(j+1)^{k}}{(k+1)^{\beta +r}k!}%
e^{\lambda \alpha (j+1)}\dint\limits_{0}^{\infty }y^{\frac{r}{\beta }}e^{y}dy
\label{26..}
\end{equation}%
Since, $\Gamma (z)=x^{z}\dint\limits_{0}^{\infty }e^{xt}t^{z-1}dt$ \ \ $,$ $%
z>0,$ $x>0,$ then%
\begin{equation}
\dint\limits_{0}^{\infty }y^{\frac{r}{\beta }}e^{y}dy=\Gamma (\frac{r}{\beta 
}+1).  \label{27.}
\end{equation}%
Substituting from (27) into (26), we get (23). This completes the proof.

\section{Maximum liklihood estimators}

In this section, we use the method of maximum likelihood to estimate the
unknown parameters of the BEMWE distribution. Consider constant values to
the parameters $\alpha ,\beta $ and $\lambda $ so, we want to estimate the
other parameters $\gamma _{1},\gamma _{2}$ and $\gamma _{3}.$ Suppose that
we have a sample of size n ,of the form $\{(x_{11},x_{21})$, $%
(x_{12},x_{22}) $,..., $(x_{1n},x_{2n})\}$ from BEMWE distribution. We use
the following notation

$I_{1}=\{x_{1i}<x_{2i}\}$, \ \ $I_{2}=\{x_{1i}>x_{2i}\}$, \ \ \ $%
I_{3}=\{x_{1i}=x_{2i}=x_{i}\}$, $\ \ I=I_{1}\cup I_{2}\cup I_{3}$,

$\left\vert I_{1}\right\vert =n_{1},$ $\ \left\vert I_{2}\right\vert =n_{2},$
$\ \left\vert I_{3}\right\vert =n_{3},$ and $n_{1}+n_{2}+n_{3}=n.$

Based on the observations, the likelihood function of the sample of size n
given by:

\begin{equation*}
l(\gamma _{1},\gamma _{2},\gamma _{3},\alpha ,\beta ,\lambda
)=\dprod\limits_{i=1}^{n_{1}}f_{1}(x_{1i},x_{2i})\dprod%
\limits_{i=1}^{n_{2}}f_{2}(x_{1i},x_{2i})\dprod%
\limits_{i=1}^{n_{3}}f_{3}(x_{i},x).
\end{equation*}

The log-likelihood function can be written as%
\begin{equation*}
L(\gamma _{1},\gamma _{2},\gamma _{3},\alpha ,\beta ,\lambda )=n_{1}\ln
\left( {\small \gamma }_{2}\left( \gamma _{1}+\gamma _{3}\right) {\small %
\lambda }^{2}{\small \beta }^{2}\right) +\dsum\limits_{i=1}^{{\small n}_{1}}(%
\frac{x_{1i}}{\alpha })^{\beta }-\lambda \alpha \dsum\limits_{i=1}^{{\small n%
}_{1}}({\small e}^{(x_{1i}/\alpha )^{{\small \beta }}}-{\small 1})\text{\ \ }
\end{equation*}%
\begin{eqnarray}
&&-\lambda \alpha \dsum\limits_{i=1}^{{\small n}_{1}}({\small e}%
^{(x_{2i}/\alpha )^{{\small \beta }}}-{\small 1})+\left( {\small \gamma }_{%
{\small 2}}{\small -1}\right) \dsum\limits_{i=1}^{{\small n}_{1}}\ln (%
{\small 1-e}^{-{\small \lambda \alpha }({\tiny e}^{(x_{2i}/\alpha )^{{\small %
\beta }}}-{\small 1})})+\left( {\small \beta -1}\right) \dsum\limits_{i=1}^{%
{\small n}_{1}}\ln (\frac{{\tiny x}_{{\tiny 1}}{\tiny x}_{{\tiny 2}}}{%
{\small \alpha }^{2}})  \notag \\
&&+\left( {\small \gamma }_{1}{\small +\gamma }_{3}{\small -1}\right)
\dsum\limits_{i=1}^{{\small n}_{1}}\ln ({\small 1-e}^{-{\small \lambda
\alpha }({\small e}^{(x_{1i}/\alpha )^{{\small \beta }}}-{\small 1}%
)})+\dsum\limits_{i=1}^{{\small n}_{1}}(\frac{x_{2i}}{\alpha })^{\beta
}+n_{2}\ln \left( {\small \gamma }_{1}\left( \gamma _{2}+\gamma _{3}\right) 
{\small \lambda }^{2}{\small \beta }^{2}\right)  \notag \\
&&+\dsum\limits_{i=1}^{{\small n}_{2}}(\frac{x_{1i}}{\alpha })^{\beta
}-\lambda \alpha \dsum\limits_{i=1}^{{\small n}_{2}}({\small e}%
^{(x_{1i}/\alpha )^{{\small \beta }}}-{\small 1})+\dsum\limits_{i=1}^{%
{\small n}_{2}}(\frac{x_{2i}}{\alpha })^{\beta }\text{\ }-\lambda \alpha
\dsum\limits_{i=1}^{{\small n}_{2}}({\tiny e}^{(x_{2i}/\alpha )^{{\small %
\beta }}}-{\small 1})  \notag \\
&&+\left( {\small \gamma }_{{\small 1}}{\small -1}\right)
\dsum\limits_{i=1}^{{\small n}_{2}}\ln ({\small 1-e}^{-{\small \lambda
\alpha }({\small e}^{(x_{1i}/\alpha )^{{\small \beta }}}-{\small 1}%
)})+\left( {\small \gamma }_{2}{\small +\gamma }_{3}{\small -1}\right)
\dsum\limits_{i=1}^{{\small n}_{2}}\ln ({\small 1-e}^{-{\small \lambda
\alpha }({\small e}^{(x_{2i}/\alpha )^{{\small \beta }}}-{\small 1})}) 
\notag \\
&&+\left( {\tiny \beta -1}\right) \dsum\limits_{i=1}^{{\small n}_{2}}\ln (%
\frac{{\tiny x}_{{\tiny 1}}{\tiny x}_{{\tiny 2}}}{{\small \alpha }^{2}}%
)+n_{3}\ln \left( \gamma _{3}{\small \lambda \beta }\right) -\lambda \alpha
\dsum\limits_{i=1}^{{\small n}_{3}}({\small e}^{(x_{i}/\alpha )^{{\small %
\beta }}}-{\small 1})+\dsum\limits_{i=1}^{{\small n}_{3}}(\frac{x_{i}}{%
\alpha })^{\beta }  \notag \\
&&+\left( {\tiny \beta -1}\right) \dsum\limits_{i=1}^{{\small n}_{3}}\ln (%
\frac{{\tiny x}_{{\tiny i}}}{{\tiny \alpha }})+\left( {\small \gamma }_{1}%
{\small +\gamma _{2}+\gamma }_{3}{\small -1}\right) \dsum\limits_{i=1}^{%
{\small n}_{3}}\ln ({\tiny 1-e}^{-{\small \lambda \alpha }({\small e}%
^{(x_{i}/\alpha )^{{\small \beta }}}-{\small 1})}).  \label{28...}
\end{eqnarray}%
Computing the first partial derivatives of (28) with respect to $\gamma
_{1},\gamma _{2}$ and $\gamma _{3}$ and setting the results equal zeros, we
get the likelihood equations as in the following form%
\begin{eqnarray}
\frac{\partial L}{\partial \gamma _{1}} &=&\frac{n_{1}}{\gamma _{1}+\gamma
_{3}}+\dsum\limits_{i=1}^{{\small n}_{1}}\ln ({\small 1-e}^{-{\small \lambda
\alpha }({\small e}^{(x_{1i}/\alpha )^{{\small \beta }}}-{\small 1})})+\frac{%
n_{2}}{\gamma _{1}}+\dsum\limits_{i=1}^{{\small n}_{2}}\ln ({\small 1-e}^{-%
{\small \lambda \alpha }({\small e}^{(x_{1i}/\alpha )^{{\small \beta }}}-%
{\small 1})})  \notag \\
&&+\dsum\limits_{i=1}^{{\small n}_{3}}\ln ({\small 1-e}^{-{\small \lambda
\alpha }({\small e}^{(x_{i}/\alpha )^{{\small \beta }}}-{\small 1})}),
\label{29..}
\end{eqnarray}%
\begin{eqnarray}
\frac{\partial L}{\partial \gamma _{2}} &=&\frac{n_{1}}{\gamma _{2}}%
+\dsum\limits_{i=1}^{{\small n}_{1}}\ln ({\small 1-e}^{-{\small \lambda
\alpha }({\small e}^{(x_{2i}/\alpha )^{{\small \beta }}}-{\small 1})})+\frac{%
n_{2}}{\gamma _{2}+\gamma _{3}}+\dsum\limits_{i=1}^{{\small n}_{2}}\ln (%
{\tiny 1-}{\small e}^{-{\small \lambda \alpha }({\small e}^{(x_{2i}/\alpha
)^{{\small \beta }}}-{\small 1})})  \notag \\
&&+\dsum\limits_{i=1}^{{\small n}_{3}}\ln ({\small 1-e}^{-{\small \lambda
\alpha }({\small e}^{(x_{i}/\alpha )^{{\small \beta }}}-{\small 1})})
\label{30.}
\end{eqnarray}%
and%
\begin{eqnarray}
\frac{\partial L}{\partial \gamma _{3}} &=&\frac{n_{1}}{\gamma _{1}+\gamma
_{3}}+\dsum\limits_{i=1}^{{\small n}_{1}}\ln ({\small 1-e}^{-{\small \lambda
\alpha }({\small e}^{(x_{1i}/\alpha )^{{\small \beta }}}-{\small 1})})+\frac{%
n_{2}}{\gamma _{2}+\gamma _{3}}+\dsum\limits_{i=1}^{{\small n}_{2}}\ln (%
{\small 1-e}^{-{\small \lambda \alpha }({\small e}^{(x_{2i}/\alpha )^{%
{\small \beta }}}-{\small 1})})  \notag \\
&&+\frac{n_{3}}{\gamma _{3}}+\dsum\limits_{i=1}^{{\small n}_{3}}\ln ({\small %
1-e}^{-{\small \lambda \alpha }({\small e}^{(x_{i}/\alpha )^{{\small \beta }%
}}-{\small 1})}).  \label{31..}
\end{eqnarray}%
To get the MLEs of the parameters $\gamma _{1},\gamma _{2}$ and $\gamma _{3}$
, we have to solve the above system of three non-linear equations. The
solution of equations (29), (30) and (31) are not easy to solve, so
numerical technique is needed to get the MLEs.

\subsection{Asymptotic confidence bounds}

In this subsection we consider the approximate confidence intervals of the
parameters $\gamma _{1},\gamma _{2}$ and $\gamma _{3}$ by using variance
covariance matrix $I_{0}^{-1}$see Lawless (2003), where $I_{0}^{-1}$ is the
inverse of the observed information matrix%
\begin{equation}
I_{0}^{-1}=-%
\begin{pmatrix}
\frac{\partial ^{2}L}{\partial \gamma _{1}^{2}} & \frac{\partial ^{2}L}{%
\partial \gamma _{1}\partial \gamma _{2}} & \frac{\partial ^{2}L}{\partial
\gamma _{1}\partial \gamma _{3}} \\ 
\frac{\partial ^{2}L}{\partial \gamma _{2}\partial \gamma _{1}} & \frac{%
\partial ^{2}L}{\partial \gamma _{2}^{2}} & \frac{\partial ^{2}L}{\partial
\gamma _{2}\partial \gamma _{3}} \\ 
\frac{\partial ^{2}L}{\partial \gamma _{3}\partial \gamma _{1}} & \frac{%
\partial ^{2}L}{\partial \gamma _{3}\partial \gamma _{3}} & \frac{\partial
^{2}L}{\partial \gamma _{3}^{2}}%
\end{pmatrix}%
^{-1}=%
\begin{pmatrix}
{\small Var(}\overset{\wedge }{{\small \gamma }_{1}}{\small )} & {\small Cov(%
\overset{\wedge }{\gamma _{1}},\gamma }_{2}{\small )} & {\small Cov(\overset{%
\wedge }{\gamma _{1}},\overset{\wedge }{\gamma _{3}})} \\ 
{\small Cov(\overset{\wedge }{\gamma _{2}},\overset{\wedge }{\gamma _{1}})}
& {\small Var(\overset{\wedge }{\gamma _{2}})} & {\small Cov(\overset{\wedge 
}{\gamma _{2}},\overset{\wedge }{\gamma _{3}})} \\ 
{\small Cov(\overset{\wedge }{\gamma _{3}},\overset{\wedge }{\gamma _{1}})}
& {\small Cov(\overset{\wedge }{\gamma _{3}},\overset{\wedge }{\gamma _{2}})}
& {\small Var(\overset{\wedge }{\gamma _{3}})}%
\end{pmatrix}%
.  \label{32}
\end{equation}%
The derivatives in $I_{0}^{-1}$ are given as follows

\begin{eqnarray*}
\frac{\partial ^{2}L}{\partial \gamma _{1}^{2}} &=&-\frac{n_{1}}{(\gamma
_{1}+\gamma _{3})^{2}}-\frac{n_{2}}{\gamma _{1}{}^{2}},\text{ \ \ \ \ \ \ }%
\frac{\partial ^{2}L}{\partial \gamma _{1}\partial \gamma _{2}}=0,\text{ \ \
\ \ \ \ }\frac{\partial ^{2}L}{\partial \gamma _{1}\partial \gamma _{3}}=-%
\frac{n_{1}}{(\gamma _{1}+\gamma _{3})^{2}}, \\
\frac{\partial ^{2}L}{\partial \gamma _{2}^{2}} &=&-\frac{n_{1}}{\gamma
_{2}^{2}}-\frac{n_{2}}{(\gamma _{2}+\gamma _{3})^{2}},\text{ \ \ \ \ \ \ \ \ 
}\frac{\partial ^{2}L}{\partial \gamma _{2}\partial \gamma _{3}}=-\frac{n_{2}%
}{(\gamma _{2}+\gamma _{3})^{2}}
\end{eqnarray*}%
and%
\begin{equation*}
\frac{\partial ^{2}L}{\partial \gamma _{3}^{2}}=-\frac{n_{1}}{(\gamma
_{1}+\gamma _{3})^{2}}-\frac{n_{2}}{(\gamma _{2}+\gamma _{3})^{2}}-\frac{%
n_{3}}{\gamma _{3}^{2}}.
\end{equation*}%
We can derive the $(1-\delta )100\%$ \ confidence intervals of the
parameters $\overset{\wedge }{{\small \gamma }_{1}},\overset{\wedge }{%
{\small \gamma }_{2}}$ and $\overset{\wedge }{{\small \gamma }_{3}}$ by
using variance covariance matrix as in the following forms%
\begin{equation*}
\overset{\wedge }{{\small \gamma }_{i}}\pm Z_{\frac{\delta }{2}}\sqrt{%
{\small Var(}\overset{\wedge }{{\small \gamma }_{i}}{\small )}}\text{ \ \ }%
,i=1,2,3.
\end{equation*}%
where $Z_{\frac{\delta }{2}}$ is the upper ($\frac{\delta }{2})$th
percentile of the standard normal distribution.

\section{Data analysis}

In this section we present the analysis of a data set and we consider a
constant value to the parameters $\alpha ,$ $\beta $ and $\lambda $\ which
take the values $0.1,$ $0.3$ and $0.05$ respectively. The data set has been
represent the American Football (National Football League) League data and
they are obtained from the matches played on three consecutive weekends in
1986. The data were first published in `Washington Post' and they are also
available in Csorgo and Welsh (1989). It is a bivariate data set, and the
variables $X_{1}$ and $X_{2}$ are as follows: $X_{1}$ represents the `game
time' to the first points scored by kicking the ball between goal posts, and 
$X_{2}$ represents the `game time' to the first points scored by moving the
ball into the end zone. These times are of interest to a casual spectator
who wants to know how long one has to wait to watch a touchdown or to a
spectator who is interested only at the beginning stages of a game.

The data (scoring times in minutes and seconds) are represented in Table 1.
Here also all the data points are divided by 100 just for computational
purposes. The variables have the following structure: (i) $X_{1}<X_{2}$
means that the first score is a field goal, (ii) $X_{1}>X_{2}$, means the
first score is an unconverted touchdown or safety, (iii) $X_{1}=X_{2}$ means
the first score is a converted touchdown.

\begin{eqnarray*}
&&\ 
\begin{tabular}{||c|c|c|c|c|c|c|c||}
\hline\hline
$X_{1}$ & $X_{2}$ & $X_{1}$ & $X_{2}$ & $X_{1}$ & $X_{2}$ & $X_{1}$ & $X_{2}$
\\ \hline\hline
2.05 & 3.98 & 8.53 & 14.57 & 2.90 & 2.90 & 1.38 & 1.38 \\ \hline
9.05 & 9.05 & 31.13 & 49.88 & 7.02 & 7.02 & 10.53 & 10.53 \\ \hline
0.85 & 0.85 & 14.58 & 20.57 & 6.42 & 6.42 & 12.13 & 12.13 \\ \hline
3.43 & 3.43 & 5.78 & 25.98 & 8.98 & 8.98 & 14.58 & 14.58 \\ \hline
7.78 & 7.78 & 13.80 & 49.75 & 10.15 & 10.15 & 11.82 & 11.82 \\ \hline
10.57 & 14.28 & 7.25 & 7.25 & 8.87 & 8.87 & 5.52 & 11.27 \\ \hline
7.05 & 7.05 & 4.25 & 4.25 & 10.40 & 10.25 & 19.65 & 10.70 \\ \hline
2.58 & 2.58 & 1.65 & 1.65 & 2.98 & 2.98 & 17.83 & 17.83 \\ \hline
7.23 & 9.68 & 6.42 & 15.08 & 3.88 & 6.43 & 10.85 & 38.07 \\ \hline
6.85 & 34.58 & 4.22 & 9.48 & 0.75 & 0.75 &  &  \\ \hline
32.45 & 42.35 & 15.53 & 15.53 & 11.63 & 17.37 &  &  \\ \hline\hline
\end{tabular}
\\
&&\text{ \ \ \ \ \ \ \textbf{Table (1).} American Football League (NFL) data}
\end{eqnarray*}

From this data, we find that the values of the unknown parameters $\overset{%
\wedge }{\gamma _{1}},$ $\overset{\wedge }{\gamma _{2}}$ and $\overset{%
\wedge }{\gamma _{3}}$ are 0.0416, 0.253 and 0.52 respectively and the
log-likelihood equals (-250.28 ). By substituting the MLE of unknown
parameters in (32), we get estimation of the variance covariance matrix as%
\begin{equation*}
I_{0}^{-1}=%
\begin{pmatrix}
0.000842 & 0.00000395 & -0.000299 \\ 
0.00000395 & 0.00394 & -0.0000939 \\ 
-0.000299 & -0.000094 & 0.00711%
\end{pmatrix}%
.
\end{equation*}%
The 95\% confidence intervals of $\overset{\wedge }{\gamma _{1}},$ $\overset{%
\wedge }{\gamma _{2}}$ and $\overset{\wedge }{\gamma _{3}}$ are (0,0.098),
(0.130,0.376) and (0.355, 0.685) respectively.\ \ \ \ \ \ \ \ \ \ \ \ \ \ \
\ \ \ \ \ \ 

\section{Conclusions}

In this paper we have introduced the bivariate exponentiated modified
Weibull extension distribution whose marginals are exponentiated modified
Weibull extension distributions. We discussed some statistical properties of
the new bivariate distribution. Maximum likelihood estimates of the new
distribution are discussed and we provided the observed Fisher information
matrix. One real data sets are analyzed using the new distribution.

\setlength{\parindent}{0in}


\begin{thebibliography}{99}
\bibitem{} Al-Khedhairi, A. and El-Gohary, A. (2008). "A new class of
bivariate Gompertz distributions" Internatinal Journal of Mathematics
Analysis, 2(5), 235 -- 253.

\bibitem{} Basu, A.P.(1971)."Bivariate failure rate". American Statistics
Association, 66,103-104.

\bibitem{} B. Gompertz, (1824) On the nature of the function expressive of
the law of human mortality and on the new mode of determining the value of
life contingencies, Philosophical Transactions of Royal Society A115, pp.
513-580.

\bibitem{} Chen, Z. (2000). "A new two-parameter lifetime distribution with
bathtub shape or increasing failure rate function". Statistics and
Probability Letters, 49, 155--161.

\bibitem{} Csorgo, S. and Welsh, A.H. (1989). "Testing for exponential and
Marshall-Olkin distribution". Journal of Statistical Planning and Inference,
vol. 23, 287-300.

\bibitem{} El-Gohary, A., Alshamrani, A. and Al-Otaibi, A. N. (2013). "The
Generalized Gompertz Distribution". Journal of Applied Mathematical
Modelling, 37(1-2), 13-24.

\bibitem{} El-Sherpieny, E. A., Ibrahim, S. A., and Bedar, R. E. (2013). "A
new bivariate generalized Gompertz distribution". Asian Journal of Applied
Sciences, 1-4, 2321 -- 0893.

\bibitem{} Johnson, N.L. and Kotz, S. (1975). "A vector valued multivariate
hazard rate". Journal of Multivariate Analysis, 5, 53-66.

\bibitem{} Kundu, D., Gupta, K. (2013)." Bayes estimation for the
Marshall--Olkin bivariate Weibull distribution". Journal of Computational
Statistics and Data Analysis, 57(1), 271--281.

\bibitem{} Kundu, D. and Gupta, R. D. (2009). "Bivariate generalized
exponential distribution". Journal of Multivariate Analysis, 100(4), 581-593.

\bibitem{} Kundu, D., Gupta, K. (2013)." Bayes estimation for the
Marshall--Olkin bivariate Weibull distribution". Journal of Computational
Statistics and Data Analysis, 57(1), 271--281.

\bibitem{} Lawless, J. F. (2003). "Statistical Models and Methods for
Lifetime Data". John Wiley and Sons, New York, 20, 1108-1113.

\bibitem{} Sarhan, A. and Balakrishnan, N. (2007). "A new class of bivariate
distributions and its mixture". Journal of the Multivariate Analysis, 98,
1508-1527.

\bibitem{} Xie, M., Tang, Y., and Goh, T. N. (2002). "A modified Weibull
extension with bathtub-shaped failure rate function". Reliability
Engineering and System Safety, 76, 279--285.
\end{thebibliography}
\end{document}